\documentclass{amsart}

\newtheorem{theorem}{Theorem}[section]

\theoremstyle{definition}

\newtheorem{example}[theorem]{Example}

\theoremstyle{remark}

\usepackage{graphicx}
\usepackage{subcaption}
\usepackage{hyperref}
\numberwithin{equation}{section}

\begin{document}

\title{Identification of 2-Bridge Links}

\author{Al\.i Sa\.it Dem\.ir}
\address{department of mathematics, \.istanbul technical university, maslak, \.istanbul 34469, turkey}
\email{demira@itu.edu.tr}
\thanks{The author was supported in part by T\"{U}B\.ITAK (The Scientific and Technological Research Council of Turkey) Grant 117F341.}






\begin{abstract}
 We find all 2-Bridge links up to 11 crossings and locate them in Thistlethwaite's link table. The splitting numbers of some links are calculated as a consequence of this identification.
\end{abstract}

\maketitle

\section{introduction}
Various properties of links are computed and listed in KnotAtlas \cite{katlas} or Knotilus \cite{knotilus} or Knotinfo \cite{knotinf} databases. The aim of this study is to identify 2-bridge links inside the list of links, so that the information about 2-bridge links can be related to their representatives in Thistlethwaite's table of links and vice versa.

\begin{figure}[h]
  \includegraphics[width=0.7\linewidth]{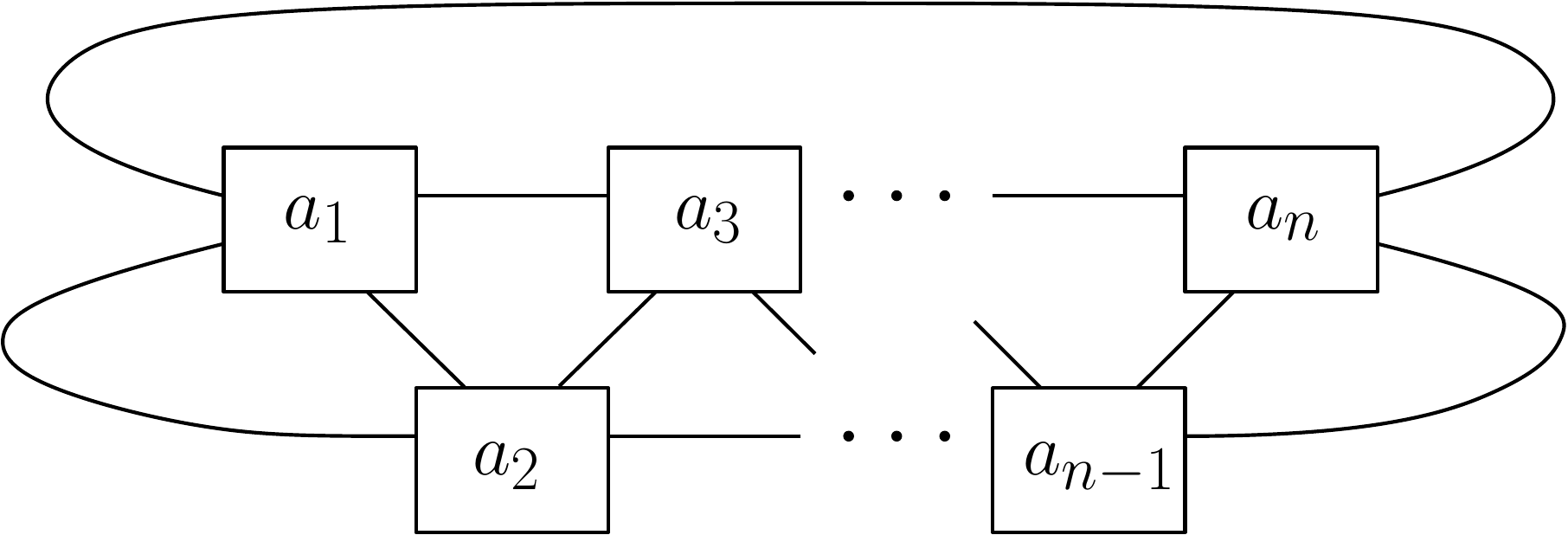}
  \caption{A 2-bridge knot or link.}
  \label{fig:2brd}
\end{figure}

Two-bridge links are links with two components which can be put into the form as in Figure \ref{fig:2brd}, where the integers $a_1,\dots,a_n$ denote the number of overcrossings (or undercrossings). 
They are also called rational links since they can be classified by the rational number $$\frac{p}{q}=a_1+\frac{1}{a_2+\frac{1}{\ddots+\frac{1}{a_n}}}$$ given by the continued fraction on the right, where $gcd(p,q)=1$. The Conway normal form $C(a)=C(a_1,\dots,a_n)$ can also be used to identify this link. Note that, for $\frac{p}{q}$ to represent a link $p$ must be even. See \cite{mur96} for more details on 2-bridge links.

 An algorithm for identifying 2-bridge links of a certain number of crossing is given in the next section. As an application, in section \ref{sec:split} we calculate the splitting numbers of some 2-bridge links that have certain type of Conway normal forms.
The following theorems, due to Schubert \cite{schu56}, are used frequently in this study to eliminate equivalent links from the list of possible combinations.
\begin{theorem}\textbf{(Schubert) (Theorem 9.3.3 of \cite{mur96})}\label{thm:sch}
The 2-bridge links $L(p,q)$ and $L(p',q')$ are equivalent as unoriented links if and only if $p=p'$ and $qq'\equiv 1\ (mod\ p)$.
\end{theorem}
\begin{theorem}\textbf{(Schubert) (Theorem 9.4.1 of \cite{mur96})}\label{thm:sch2}
If the orientation of one component of the 2-bridge link $L(p,q)$, where both $p>0$ and $q>0$, is reversed the resulting link is equivalent to $L(p,q-p)=L^*(p,p-q)$.
\end{theorem}
In Example \ref{ex1} below, we exhibit how these theorems are used to find equivalent (up to orientation) links of 7 crossings.

\section{Identifying and listing the links}

In this section an algorithm to find all 2-bridge links with $n$ crossings is outlined and the results for links up to 11 crossing are listed in section \ref{sec:lists}. These links are matched with their Thistlethwaite's Id to make other data about these links relatable to their Conway normal forms or rational representations. For a similar study on the 2-bridge knots, see De Wit's paper\cite{DeWit}.
In order to find and identify 2-bridge links up to a certain number of crossings $n$:\\
\begin{enumerate}
	\item Find all permutations of positive integers that add up to $n$,
	\item Calculate the rational number $\frac{p}{q}$ for each such permutation,
	\item Rule out the permutations giving knots instead of links, by checking the parity of $p$,
	\item Shorten the list by picking one representative of permutations with equal continued fraction or equivalent link according to Schubert's criteria (Theorem \ref{thm:sch} and \ref{thm:sch2}),
	\item Identify the Gauss code of the link by looking at its Conway form and locate it in Thistlethwaite's link table.
\end{enumerate}

 The Thistlethwaite's link table lists links up to orientations and mirror images. This means in the tables below a link may represent 2 (if there is only one orientation which is the case the link has a palindromic Conway form) or 4 (two orientations of both the link and its mirror) links.

\begin{table}[h]

\begin{tabular}{ |p{1cm}|p{1cm}|p{2cm}|p{0.05cm}|p{1cm}|p{1cm}|p{2cm}|  }
 
 \hline
 ID&Link $(p,q)$&Conway Form&& ID&Link $(p,q)$&Conway Form \\
 \hline
$E_{1 }$	 & $(	14	,	3	)$	& $[	4	,	1	,	1	,	1					]$&	&$E_{ 11}$	 & $(	16	,	9	)$	& $[	1	,	1	,	3	,	1	,	1			]$	\\	\hline
$E_{ 2}$	 & $(	14	,	3	)$	& $[	4	,	1	,	2							]$	&&$E_{ 12}$	 & $(	16	,	9	)$	& $[	1	,	1	,	3	,	2					]$	\\	\hline
$E_{ 3}$	 & $(	14	,	5	)$	& $[	2	,	1	,	3	,	1					]$	&&$E_{ 13}$	 & $(	18	,	5	)$	& $[	3	,	1	,	1	,	1	,	1			]$	\\	\hline
$E_{ 4}$	 & $(	14	,	5	)$	& $[	2	,	1	,	4							]$	&&$E_{ 14}$	 & $(	18	,	5	)$	& $[	3	,	1	,	1	,	2					]$	\\	\hline
$E_{ 5}$	 & $(	14	,	9	)$	& $[	1	,	1	,	1	,	3	,	1			]$	&&$E_{ 15}$	 & $(	18	,	7	)$	& $[	2	,	1	,	1	,	2	,	1			]$	\\	\hline
$E_{ 6}$	 & $(	14	,	9	)$	& $[	1	,	1	,	1	,	4					]$	&&$E_{ 16}$	 & $(	18	,	7	)$	& $[	2	,	1	,	1	,	3					]$	\\	\hline
$E_{ 7}$	 & $(	14	,	11	)$	& $[	1	,	3	,	1	,	1	,	1	  ]$	&&$E_{ 17}$	 & $(	18	,	11	)$	& $[	1	,	1	,	1	,	1	,	2	,	1	]$	\\	\hline
$E_{ 8}$	 & $(	14	,	11	)$	& $[	1	,	3	,	1	,	2				]$	&&$E_{ 18}$	 & $(	18	,	11	)$	& $[	1	,	1	,	1	,	1	,	3			]$	\\	\hline
$E_{ 9}$	 & $(	16	,	7	)$	& $[	2	,	3	,	1	,	1					]$	&&$E_{ 19}$	 & $(	18	,	13	)$	& $[	1	,	2	,	1	,	1	,	1	,	1	]$	\\	\hline
$E_{10 }$	 & $(	16	,	7	)$	& $[	2	,	3	,	2							]$	&&$E_{ 20}$	 & $(	18	,	13	)$	& $[	1	,	2	,	1	,	1	,	2			]$	\\	\hline
\end{tabular}
\caption{All possible combinations of 7 crossings}
\label{tab:tab}
\end{table}

\begin{example}\label{ex1}
In Table \ref{tab:tab} we list all 20 possible 2-bridge links with 7 crossings and explain why these are represented by only 3 links in Table \ref{tab:tab7}. Consider the links $E_1$ to $E_8$. The links with continued fractions leading to equal rational number are equivalent, hence we have $E_1=E_2$, $E_3=E_4$, $E_5=E_6$ and $E_7=E_8$. By Theorem \ref{thm:sch} $E_1=E_3$ since $3\cdot5\equiv 1\ (mod \ 14)$. Similarly $E_5=E_7$. These equivalences assume that the links are  not oriented. Suppose $E_3$ is given the standard orientation. According to Theorem \ref{thm:sch2}, if the orientation of one of the components of $E_3$ is reversed the resulting link is equivalent (as oriented links) to $E_5$. Since Thistlethwaite's link table lists links up to orientations all these 8 links are represented by one link, namely $L_7A_6$. Similarly the links $E_9$ to $E_{12}$ and $E_{13}$ to $E_{20}$ are represented by $L_7A_4$ and $L_7A_5$, respectively.
\end{example}
\section{splitting numbers}\label{sec:split}
Although the classification of 2-bridge links is complete, various characteristics or local properties of these links are actively studied. In this section we will mention one of such invariants, namely the splitting number of a link.

The splitting number $sp(L)$ of the link $L$ is defined to be the minimum number of crossing changes between different components of $L$ to convert $L$ into a split link. For more information about splitting numbers see \cite{ccz16} and \cite{cfp17} and references therein. In \cite{cfp17} authors calculate the splitting number of links up to 9 crossings. They use 5 methods based on covering properties or Alexander invariants, case by case, for determining the splitting numbers. The other study \cite{ccz16}, due to Cimasoni et. al., calculates the splitting number by looking at the signature and nullity of the links. As a consequence of their main result the following theorem is about the splitting number of certain 2-bridge links:

\begin{theorem}\textbf{(Theorem 4.7 of \cite{ccz16})}
The splitting number of the 2-bridge link $C(2a_1,b_1,\dots,2a_{n-1},b_{n-1},2a_n)$ is $a_1+a_2+\dots+a_n$, where all $a_i$ and $b_i$ are positive integers.
\label{thm4.7}
\end{theorem}
\begin{table}[h]
\begin{tabular}{ |p{1.5cm}|p{1cm}|p{0.05cm}|p{1.5cm}|p{1cm}| p{0.05cm}|p{1.5cm}|p{1cm}|  }
 
 \hline
 Link & $sp(L)$ & & Link &$sp(L)$& & Link &$sp(L)$ \\ \hline
$L_{10}A_{48}$&2&&$L_{11}A_{132}$& 2&&$L_{11}A_{299}$& 4\\ \hline
$L_{10}A_{64}$&3&&$L_{11}A_{194}$& 3&&$L_{11}A_{312}$& 4\\ \hline
$L_{10}A_{75}$&3&&$L_{11}A_{206}$& 3&&$L_{11}A_{319}$& 4\\ \hline
$L_{10}A_{87}$&3&&$L_{11}A_{222}$& 3&&$L_{11}A_{360}$& 5\\ \hline
$L_{10}A_{89}$&4&&$L_{11}A_{263}$& 4&&$L_{11}A_{372}$& 5\\ \hline
$L_{10}A_{98}$&4&&$L_{11}A_{278}$& 4&&& \\ \hline
$L_{10}A_{102}$&4&&$L_{11}A_{289}$& 4&&& \\ \hline
\end{tabular}
\caption{Splitting numbers of some 2-bridge links}
\label{tab:split}
\end{table}

 Our calculations below reveal which alternating two component links are 2-bridge links. Therefore one can also calculate the splitting numbers of the following links, which turn out to be 2-bridge links with the desired Conway form, using Theorem \ref{thm4.7}: 
$$L_5A_1,\ L_6A_1,\ L_7A_4,\ L_7A_6,\ L_8A_6,\ L_8A_8,\ L_8A_{11},\ L_9A_{18},\ L_9A_{26},\ L_9A_{30},\ L_9A_{36},\ L_9A_{40}$$
Besides the above links that already appear with splitting numbers in \cite{cfp17}, we calculate the splitting numbers of those links in table \ref{tab:split} with 10 and 11 crossings that fit into the Conway form of Theorem \ref{thm4.7}:

We note that the splitting number of $L_{11}A_{372}$ was also calculated in \cite{ccz16} (Example 4.4), but not making use of the fact that Theorem \ref{thm4.7} also applies to this link.
\section*{acknowledgments}
The author would like to express his gratitude to Alexander Degtyarev for introducing the subject and commenting on an earlier draft of the paper.

\section{Two-bridge links up to 11 crossings}\label{sec:lists}

\begin{table}[h]
\begin{tabular}{ |p{3cm}|p{3cm}|p{3cm}|  }
 
 \hline
 The Link $\frac{p}{q}$ or $(p,q)$&Conway Form & Thistlethwaite's Id\\
 \hline
(4,1)  & [4] or [3,1] &$L_4A_1$\\ \hline
 (8,3)  & [2,1,2]   &$L_5A_1$\\ \hline
\end{tabular}
\caption{2-Bridge Links of 4 and 5 Crossings}

\end{table}

\begin{table}[h]
\begin{tabular}{ |p{3cm}|p{3cm}|p{3cm}|  }
 
 \hline
 The Link $\frac{p}{q}$ or $(p,q)$&Conway Form & Thistlethwaite's Id\\
 \hline
(6,1)  & [6] or [5,1]  &$L_6A_3$\\
 
 \hline
(10,3)  & [3,3]   &$L_6A_2$\\

\hline
(12,5)  & [2,2,2]   &$L_6A_1$\\ 
\hline
\end{tabular}
\caption{2-Bridge Links of 6 Crossings}

\end{table}

\begin{table}[h]

\begin{tabular}{ |p{3cm}|p{3cm}|p{3cm}|  }
 
 \hline
 The Link $\frac{p}{q}$ or $(p,q)$&Conway Form & Thistlethwaite's Id\\
 \hline
(14,5)  & [2,1,4]   &$L_7A_6$\\ \hline

  (16,7)  & [2,3,2]   &$L_7A_4$\\ \hline
(18,5)  & [3,1,1,2]   &$L_7A_5$\\
 
 \hline
\end{tabular}
\caption{2-Bridge Links of 7 Crossings}
\label{tab:tab7}
\end{table}

\begin{table}[h]
\begin{tabular}{ |p{3cm}|p{3cm}|p{3cm}|  }
 
 \hline
 The Link $\frac{p}{q}$ or $(p,q)$&Conway Form & Thistlethwaite's Id\\
 \hline
 (8,1)  & [8] or [7,1] &$L_8A_{14}$\\ 
 \hline
  (16,3)  & [5,2,1]   &$L_8A_{12}$\\
 \hline
(20,9)  & [2,4,2]   &$L_8A_6$\\ 
 \hline
(22,5)  & [4,2,2]   &$L_8A_{11}$\\ 
 \hline
  (24,7)  & [3,2,3]   &$L_8A_{13}$\\
 \hline
(26,7)  & [3,1,2,2]   &$L_8A_{10}$\\ 
 \hline
(30,11)  & [2,1,2,1,2]   &$L_8A_8$\\ 
 \hline
  (34,13)  & [2,1,1,1,1,2]   &$L_8A_9$\\
 \hline
\end{tabular}
\caption{2-Bridge Links of 8 Crossings}

\end{table}

\begin{table}[h]
\begin{tabular}{ |p{3cm}|p{3cm}|p{3cm}|  }
 
 \hline
 The Link $\frac{p}{q}$ or $(p,q)$&Conway Form & Thistlethwaite's Id\\
 \hline
 (20,7)  & [2,1,6]   &$L_9A_{36}$\\ 
 \hline
(24,5)  & [4,1,4]   &$L_9A_{40}$\\ 
 \hline
(24,11)  & [2,5,2]   &$L_9A_{18}$\\ 
 \hline
(28,11)  & [2,1,1,5]   &$L_9A_{39}$\\ 
 \hline
(30,7)  & [4,3,2]   &$L_9A_{30}$\\ 
\hline
(32,7)  & [4,1,1,3]   &$L_9A_{38}$\\ 
 \hline
(34,9)  & [3,1,3,2]   &$L_9A_{25}$\\ 
 \hline
(36,11)  & [3,3,1,2]   &$L_9A_{34}$\\ 
 \hline
(40,11)  & [3,1,1,1,3]   &$L_9A_{35}$\\ 
 \hline
(44,13)  & [3,2,1,1,2]   &$L_9A_{37}$\\ 
 \hline
(46,17)  & [2,1,2,2,2]   &$L_9A_{26}$\\ 
 \hline
(50,19)  & [2,1,1,1,2,2]   &$L_9A_{27}$\\ 
 \hline
\end{tabular}
\caption{2-Bridge Links of 9 Crossings}
\end{table}

\begin{table}[h]
\begin{tabular}{ |p{3cm}|p{3cm}|p{3cm}|  }
 
 \hline
 The Link $\frac{p}{q}$ or $(p,q)$&Conway Form & Thistlethwaite's Id\\
 \hline
(10,1)  & [10] or [9,1] &$L_{10}A_{118}$ \\  \hline
(22,3)  & [7,3]   &$L_{10}A_{114}$ \\  \hline
(26,5)  & [5,5]   &$L_{10}A_{120}$ \\  \hline
(28,13) & [2,6,2]   &$L_{10}A_{48}$\\  \hline
(32,5)  & [6,2,2]   &$L_{10}A_{98}$ \\  \hline
(38,9)  & [4,4,2]   &$L_{10}A_{75}$ \\  \hline
(40,7)  & [5,1,2,2]   &$L_{10}A_{97}$ \\  \hline
(40,9)  & [4,2,4]   &$L_{10}A_{102}$ \\  \hline
(42,11)  & [3,1,4,2]   &$L_{10}A_{73}$ \\  \hline
(42,13)  & [3,4,3]   &$L_{10}A_{115}$ \\  \hline
(48,11)  & [4,2,1,3]   &$L_{10}A_{100}$ \\  \hline
(48,17)  & [2,1,4,1,2]   &$L_{10}A_{89}$ \\  \hline
(52,11)  & [4,1,2,1,2]   &$L_{10}A_{99}$ \\  \hline
(56,15)  & [3,1,2,1,3]   &$L_{10}A_{101}$ \\  \hline
(56,17)  & [3,3,2,2]   &$L_{10}A_{94}$ \\  \hline
(58,17)  & [3,2,2,3]   &$L_{10}A_{116}$ \\  \hline
(60,13)  & [4,1,1,1,1,2]   &$L_{10}A_{93}$ \\  \hline
(62,23)  & [2,1,2,3,2]   &$L_{10}A_{64}$ \\  \hline
(64,19)  & [3,2,1,2,2]   &$L_{10}A_{96}$ \\  \hline
(64,23)  & [2,1,3,1,1,2]   &$L_{10}A_{90}$ \\  \hline
(66,25)  & [2,1,1,1,3,2]   &$L_{10}A_{65}$ \\  \hline
(68,19)  & [3,1,1,2,1,2]   &$L_{10}A_{92}$ \\  \hline
(70,29)  & [2,2,2,2,2]   &$L_{10}A_{87}$ \\  \hline
(74,31)  & [2,2,1,1,2,2]   &$L_{10}A_{83}$ \\  \hline
(76,21)  & [3,1,1,1,1,1,2]   &$L_{10}A_{88}$ \\  \hline
(80,31)  & [2,1,1,2,1,1,2]   &$L_{10}A_{91}$ \\  \hline
 
\end{tabular}
\caption{2-Bridge Links of 10 Crossings}
\end{table}

\begin{table}[h]
\begin{tabular}{ |p{3cm}|p{3cm}|p{3cm}|  }
 
 \hline
 The Link $\frac{p}{q}$ or $(p,q)$&Conway Form & Thistlethwaite's Id\\
 \hline
$(	26	,	3	)$	&	$[	8	,	1	,	2															]$	&	$L_{11}A_{	360	}$	 \\  \hline
$(	32	,	15	)$	&	$[	2	,	7	,	2															]$	&	$L_{11}A_{	132	}$	 \\  \hline
$(	34	,	5	)$	&	$[	6	,	1	,	4															]$	&	$L_{11}A_{	372	}$	 \\  \hline
$(	38	,	5	)$	&	$[	7	,	1	,	1		,	2	,											]$	&	$L_{11}A_{	367	}$	 \\  \hline
$(	44	,	7	)$	&	$[	6	,	3	,	2															]$	&	$L_{11}A_{	278	}$	 \\  \hline
$(	46	,	7	)$	&	$[	6	,	1	,	1		,	3	,											]$	&	$L_{11}A_{	364	}$	 \\  \hline
$(	46	,	11	)$	&	$[	4	,	5	,	2		,													]$	&	$L_{11}A_{	206	}$	 \\  \hline
$(	50	,	9	)$	&	$[	5	,	1	,	1		,	4	,											]$	&	$L_{11}A_{	365	}$	 \\  \hline
$(	50	,	13	)$	&	$[	3	,	1	,	5		,	2	,											]$	&	$L_{11}A_{	192	}$	 \\  \hline
$(	52	,	9	)$	&	$[	5	,	1	,	3		,	2	,											]$	&	$L_{11}A_{	275	}$	 \\  \hline
$(	54	,	17	)$	&	$[	3	,	5	,	1		,	2	,											]$	&	$L_{11}A_{	355	}$	 \\  \hline
$(	56	,	13	)$	&	$[	4	,	3	,	4															]$	&	$L_{11}A_{	319	}$	 \\  \hline
$(	58	,	11	)$	&	$[	5	,	3	,	1		,	2												]$	&	$L_{11}A_{	371	}$	 \\  \hline
$(	62	,	11	)$	&	$[	5	,	1	,	1		,	1	,	3										]$	&	$L_{11}A_{	358	}$	 \\  \hline
$(	62	,	13	)$	&	$[	4	,	1	,	3		,	3												]$	&	$L_{11}A_{	369	}$	 \\  \hline
$(	64	,	15	)$	&	$[	4	,	3	,	1		,	3												]$	&	$L_{11}A_{	302	}$	 \\  \hline
$(	70	,	13	)$	&	$[	5	,	2	,	1		,	1	,	2										]$	&	$L_{11}A_{	362	}$	 \\  \hline
$(	72	,	19	)$	&	$[	3	,	1	,	3		,	1	,	3										]$	&	$L_{11}A_{	298	}$	 \\  \hline
$(	74	,	23	)$	&	$[	3	,	4	,	1		,	1	,	2										]$	&	$L_{11}A_{	368	}$	 \\  \hline
$(	76	,	23	)$	&	$[	3	,	3	,	3		,	2												]$	&	$L_{11}A_{	260	}$	 \\  \hline
$(	76	,	27	)$	&	$[	2	,	1	,	4		,	2	,	2										]$	&	$L_{11}A_{	263	}$	 \\  \hline
$(	78	,	17	)$	&	$[	4	,	1	,	1		,	2	,	3										]$	&	$L_{11}A_{	366	}$	 \\  \hline
$(	78	,	29	)$	&	$[	2	,	1	,	2		,	4	,	2										]$	&	$L_{11}A_{	194	}$	 \\  \hline
$(	80	,	17	)$	&	$[	4	,	1	,	2		,	2	,	2										]$	&	$L_{11}A_{	312	}$	 \\  \hline
$(	82	,	23	)$	&	$[	3	,	1	,	1		,	3	,	3										]$	&	$L_{11}A_{	356	}$	 \\  \hline
$(	82	,	31	)$	&	$[	2	,	1	,	1		,	1	,	4	,	2								]$	&	$L_{11}A_{	196	}$	 \\  \hline
$(	84	,	19	)$	&	$[	4	,	2	,	2		,	1	,	2										]$	&	$L_{11}A_{	299	}$	 \\  \hline
$(	84	,	25	)$	&	$[	3	,	2	,	1		,	3	,	2										]$	&	$L_{11}A_{	271	}$	 \\  \hline
$(	86	,	25	)$	&	$[	3	,	2	,	3		,	1	,	2										]$	&	$L_{11}A_{	361	}$	 \\  \hline
$(	88	,	19	)$	&	$[	4	,	1	,	1		,	1	,	2	,	2								]$	&	$L_{11}A_{	280	}$	 \\  \hline
$(	92	,	21	)$	&	$[	4	,	2	,	1		,	1	,	1	,	2								]$	&	$L_{11}A_{	305	}$	 \\  \hline
$(	92	,	33	)$	&	$[	2	,	1	,	3		,	1	,	2	,	2								]$	&	$L_{11}A_{	264	}$	 \\  \hline
$(	94	,	39	)$	&	$[	2	,	2	,	2		,	3	,	2										]$	&	$L_{11}A_{	222	}$	 \\  \hline
$(	98	,	27	)$	&	$[	3	,	1	,	1		,	1	,	2	,	3								]$	&	$L_{11}A_{	359	}$	 \\  \hline
$(	98	,	41	)$	&	$[	2	,	2	,	1		,	1	,	3	,	2								]$	&	$L_{11}A_{	221	}$	 \\  \hline
$(	100	,	27	)$	&	$[	3	,	1	,	2		,	2	,	1	,	2								]$	&	$L_{11}A_{	297	}$	 \\  \hline
$(	100	,	39	)$	&	$[	2	,	1	,	1		,	3	,	2	,	2								]$	&	$L_{11}A_{	284	}$	 \\  \hline
$(	104	,	29	)$	&	$[	3	,	1	,	1		,	2	,	2	,	2								]$	&	$L_{11}A_{	272	}$	 \\  \hline
$(	106	,	31	)$	&	$[	3	,	2	,	2		,	1	,	1	,	2								]$	&	$L_{11}A_{	363	}$	 \\  \hline
$(	108	,	29	)$	&	$[	3	,	1	,	2		,	1	,	1	,	1	,	2						]$	&	$L_{11}A_{	300	}$	 \\  \hline
$(	112	,	31	)$	&	$[	3	,	1	,	1		,	1	,	1	,	2	,	2						]$	&	$L_{11}A_{	262	}$	 \\  \hline
$(	112	,	41	)$	&	$[	2	,	1	,	2		,	1	,	2	,	1	,	2						]$	&	$L_{11}A_{	289	}$	 \\  \hline
$(	116	,	45	)$	&	$[	2	,	1	,	1		,	2	,	1	,	2	,	2						]$	&	$L_{11}A_{	266	}$	 \\  \hline
$(	128	,	47	)$	&	$[	2	,	1	,	2		,	1	,	1	,	1	,	1	,	2				]$	&	$L_{11}A_{	247	}$	 \\  \hline
$(	144	,	55	)$	&	$[	2	,	1	,	1		,	1	,	1	,	1	,	1	,	1	,	2		]$	&	$L_{11}A_{	248	}$	 \\  \hline

\end{tabular}
\caption{2-Bridge Links of 11 Crossings}
\end{table}
\clearpage
\bibliographystyle{amsplain}

\end{document}